\documentclass[xits,reqno]{amsart}
\usepackage[active]{srcltx}
\usepackage[colorlinks=true,citecolor=black,linkcolor=black,urlcolor=blue]{hyperref}
\def\qed{\hfill \rule{4pt}{7pt}}
\def\pf{\noindent {\it Proof.} }
\newcommand{\poq}[2]{(#1;q)_{#2}}

\def\qed{\hfill \rule{4pt}{7pt}}
\def\pf{\noindent {\it Proof.} }

\usepackage{latexsym}
\usepackage{amsmath,amsthm}
\usepackage{graphicx}
\usepackage{amssymb,amsmath}
\newtheorem{dl}{Theorem}[section]
\newtheorem{lz}[dl]{Example}

\newtheorem{tl}[dl]{Corollary}
\newtheorem{yl}[dl]{Lemma}
\newtheorem{dy}[dl]{Definition}
\newtheorem{rk}[dl]{Remark}
\newtheorem{xinzhi}[dl]{Proposition}
\newtheorem{cx}[dl]{Conjecture}

\newenvironment{sourcelist}
  {\par\bgroup
   \hbox to\hsize{\hss\vrule width 0.8pt height 0cm depth 0.25cm%
                      \vrule width 1.025\hsize height 0cm depth 0.8pt%
                      \vrule width 0.8pt height 0cm depth 0.25cm\hss}
   \nopagebreak
   \parindent=0pt\small
   \tt\obeyspaces\obeylines\nopagebreak}
  {\nopagebreak\par\nopagebreak
   \hbox to\hsize{\hss\vrule width 0.8pt height 0.25cm depth 0cm%
                      \vrule width 1.025\hsize height 0.8pt depth 0cm%
                      \vrule width 0.8pt height 0.25cm depth 0cm\hss}
   \egroup\medskip}

\begin{document}
\title{Proof of the $(\alpha,\beta)$-inversion formula\\ conjectured by Hsu and Ma}
\dedicatory{\textsc{Jin Wang}~$^{a,1}$
\ and \ \textsc{Xinrong Ma}~$^{b,2},$\\
         $^a$~
Department of Mathematics, Zhejiang Normal  University\\
Jinhua 321004,~P.~R.~China\\[1mm]
$^b$~
Department of Mathematics, Soochow University\\
Suzhou 215006,~P.~R.~China}
\thanks{E-mail addresses:$^1$~\emph{jinwang@zjnu.edu.cn}
        and~$^2$~\emph{xrma@suda.edu.cn}}
\thanks{$^2$~Corresponding author. This work was supported by NSF of Zhejiang Province (Grant~No.~LQ20A010004) and by NSF of China (Grant~No.~11971341 and 12001492)}


\begin{abstract}In light of the well-known fact that the $n$th divided difference of
any polynomial of  degree $m$ must be zero while $m<n$,
the present paper proves the $(\alpha,\beta)$-inversion formula conjectured by Hsu and Ma [J. Math. Res. $\&$ Exposition  25(4) (2005) 624]. As applications of $(\alpha,\beta)$-inversion, we not only recover some known matrix inversions due to  Gasper, Schlosser, and Warnaar, but also fin  three  new matrix inversions related to elliptic divisibility sequence and theta functions.
\\

{\sl This paper is dedicated to the memory of Professor L. C. Hsu}
\end{abstract}

\maketitle
\markboth{J. Wang and X. Ma}{Proof of the $(\alpha,\beta)$-inversion formula conjectured by Hsu and Ma}

\parskip 7pt

\section{Introduction}
Throughout this paper, all operations are carried out on the complex field $\mathbb{C}$.  Recall that
    $F=(F(n,k))_{n,k \in {\Bbb Z}}$ is  an   infinite-dimensional
    lower-triangular matrix over $\mathbb{C}$, often denoted by $(F(n,k))_{n\geq k \in {\Bbb Z}}$, provided that each entry $F(n,k)=0$ unless $\ n\geq k.$ The matrix
$G=(G(n,k))_{n, k \in {\Bbb Z}}$
 is the inverse matrix of $F$  if
 \begin{align}
 \displaystyle\sum_{k\leq i\leq n}F(n,i)G(i,k)=\delta_{n,k}\,\,\mbox{for all}\ n, k\in {\Bbb Z}, \label{matrixinvdef}
 \end{align}
 where $\delta_{n,k}$ denotes the usual Kronecker delta, $\mathbb{Z}$  denotes the
set of integers.  A pair of such matrices, as pointed out by Henrici \cite{Henrici} and Gessel and Stanton
 \cite[p.175,\S 2]{111} independently,  is equivalent to the Lagrange inversion
formula  and  is often called  an inversion
formula or  a reciprocal relation in the context of Combinatorics.
In what follows,
 we call such a pair of matrices $F$ and $G$ with the reciprocal relation  a {\sl matrix inversion}.
 As many facts have shown that matrix inversions, called the inverse technique by  Chu and Hsu, play
very important roles in deriving  summation and transformation
formulas of various hypergeometric series. The reader may consult
\cite{chu,chu1,chu2, 55,101,111,0020,005,milne,1000} for more details.

It is worth noting that in the \cite{0020} Ma established
\begin{dl}[The $(f,g)$-inversion formula]\label{fg}
Preserve the above notation and assumptions. Suppose further $g(x,y)$ is anti-symmetric, i.e., $g(x,y)=-g(y,x)$.
 Let $F=(F(n,k))_{n\geq k\in {\Bbb Z}}$ and $G=(G(n,k))_{n\geq k\in {\Bbb Z}}$ be two matrices
with entries given by
\begin{subequations}
\begin{align}
F(n,k)&=\frac{\prod_{i=k}^{n-1}f(x_i,b_k)}
{\prod_{i=k+1}^{n}g(b_i,b_k)}\qquad\mbox{and}\\
G(n,k)&=
\frac{f(x_k,b_k)}{f(x_n,b_n)}\frac{\prod_{i=k+1}^{n}f(x_i,b_n)}
{\prod_{i=k}^{n-1}g(b_i,b_n)},\quad\mbox{respectively}.
\end{align}
\end{subequations}
 Then $F=(F(n,k))_{n\geq k\in {\Bbb Z}}$ and $G=(G(n,k))_{n\geq k\in {\Bbb Z}}$ is a
matrix inversion if and only if for all $a,b,c,x\in \mathbb{C}$,
there holds
\begin{align}
f(x,a)g(b,c)+f(x,b)g(c,a)+f(x,c)g(a,b)=0.
\end{align}
\end{dl}

As it turns out, the $(f,g)$-inversion formula provides a general framkwork for many existing matrix inversions. Shortly afterward, with a motivation to extend  the valid range of  the $(f,g)$-inversion formula to arbitrary discrete sequences, Hsu and Ma  \cite{hsuma}  proposed  a discrete analogue of Theorem \ref{fg} and inquired for any quick proof.  However,  it remains unproved until now for lack of the arbitrariness of continuous variables used in the $(f,g)$-inversion formula.
\begin{cx}[The $(\alpha,\beta)$-inversion formula: Hsu and Ma \cite{hsuma}]  \label{Thmainnewab}
Let $\{\alpha_{n,k}\}_{n,k\in \Bbb{Z}}$ and $\{\beta_{n,k}\}_{n,k\in \Bbb{Z}}$ be two arbitrary double
index  sequences over $\mathbb{C}$ such that none of the terms
$\alpha_{n,n}$ or $\beta_{n,k}$ is zero, and $\beta_{n,k}$ is
antisymmetric, i.e., $\beta_{n,k}=-\beta_{k,n}$. Let $F=(F(n,k))_{n\geq k\in {\Bbb Z}}$ and $G=(G(n,k))_{n\geq k\in {\Bbb Z}}$  be two infinite-dimensional
    lower-triangular matrices with entries given by
 \begin{subequations}\label{quintuplenew}
\begin{align}
&F(n,k)=\frac{\prod_{i=k}^{n-1} \alpha_{i,k}}
{\prod_{i=k+1}^{n}\beta_{i,k}}\label{cc12}\qquad\mbox{and}\\
&G(n,k)=\frac{\alpha_{k,k}}{\alpha_{n,n}}
\frac{\prod_{i=k+1}^{n} \alpha_{i,n}}
{\prod_{i=k}^{n-1}\beta_{i,n}}.\label{dd12}
\end{align}
\end{subequations} Then $F=(F(n,k))_{n\geq k\in {\Bbb Z}}$ and $G=(G(n,k))_{n\geq k\in {\Bbb Z}}$
is a matrix inversion if and only if
 for arbitrary integers $n,k,l,m,$ there holds
\begin{align}
\alpha_{n,l}\beta_{m,k}+ \alpha_{n,m}\beta_{k,l}+
 \alpha_{n,k}\beta_{l,m}=0.\label{cond20}
\end{align}
\end{cx}

In what follows, we refer to this conjecture as  the $(\alpha,\beta)$-inversion formula.  The theme of this paper is to show

\begin{dl}\label{1.33} \eqref{cond20} is sufficient but not necessary for Conjecture \ref{Thmainnewab}.
\end{dl}

 Our argument mainly relies on the following   general matrix
inversion.

\begin{yl}\label{1.3} Let $\{a_n\}_{n\in {\Bbb Z}},\{b_n\}_{n\in {\Bbb Z}}$ and $\{s_n\}_{n\in {\Bbb Z}},\{m_n\}_{n\in {\Bbb Z}}$ be four arbitrary
sequences over $\mathbb{C}$ such that none of the terms both $a_n$
and $b_n$ is zero, $s_n$ are distinct from each other. Let $F=(F(n,k))_{n\geq k\in {\Bbb Z}}$ and
$G=(G(n,k))_{n\geq k\in {\Bbb Z}}$ be two infinite-dimensional
    lower-triangular matrices with entries given by
respectively\begin{subequations}
 \begin{align}
F(n,k)=\frac{b_n}{b_k} \prod_{i=k+1}^{n} \frac{m_i(
s_k-s_{i-1}+a_{i-1}b_{i-1}\,m_{i-1}) }
      { s_k - s_i  }
\end{align}
and
\begin{align}
G(n,k)=\frac{a_{k}}{a_n}\prod_{i = k}^{n-1}\frac{\,m_i\,
      (s_n-s_{i+1}+
        a_{i+1}\,b_{i+1}\,m_{i+1})}{
      s_n-s_i}.
\end{align}
\end{subequations} Then $F$ and $G$ is a matrix inversion.
\end{yl}

Several notation on convention are needed. Hereafter,
 any product of the form $\prod_{i=k}^na_i$ for  $k,n\in {\Bbb Z}$ is defined by
(cf.\cite{10})
\begin{align}
\prod_{i=k}^na_i:= \left\{\begin{array}{ll} a_{k}a_{k+1}\cdots
a_{n},& n\geq k;\\
1,& n=k-1;\\
 1/(a_{n+1}a_{n+2}\cdots a_{k-1}),& n\leq k-2.
\end{array}\label{convention}
\right.
\end{align}
As for $q$-series, we employ the following standard notations for
the $q$-shifted factorials: for any integers $m\geq 1, n\in {\Bbb Z}$,
\begin{align*}
 &  (a;q)_\infty:=\prod_{i=0}^{\infty}(1-aq^i),\,\,(a;q)_n:=\frac{(a;q)_\infty}{(aq^n;q)_\infty};\\
 & (a_1,a_2,\ldots,a_m;q)_n:=(a_1;q)_n(a_2;q)_n\cdots(a_m;q)_n.
\end{align*}
As for theta and elliptic hypergeometric series, we adopt the  standard concepts from \cite[p.304, (11.2.5)/(11.2.6)]{10} for  Jacobi's theta function and  the theta analogue of
the $q$-shifted factorial, as follows:
\begin{align*}
   \theta(x;q):=(x,q/x;q)_\infty,\quad (x;q,p)_n:=\prod_{k=0}^{n-1}\theta(xq^k;p)
            \end{align*}
            as well as their multivariate analogues
\begin{align*}
       \theta(a_1,a_2,\ldots,a_m;q)&:=\theta(a_1;q)\theta(a_2;q)\cdots \theta(a_m;q);\\
       (a_1,a_2,\ldots,a_m;q,p)_n&:=(a_1;q,p)_n(a_2;q,p)_n\cdots (a_m;q,p)_n.
       \end{align*}
Our paper is organized as follows.   Section 2 is devoted to the proof of  Lemma \ref{1.3}. It is based on the well-known fact that the $n$-th divided difference of
any polynomial of degree $m$ must be zero while $m <n$. In the section 3,  we introduce  the so-called triple sum identity and the quintuple sum identity and show they are equivalent to each others. By using the equivalency of these two identities and their relationship with Lemma  \ref{1.3},  we finally achieve the proof of Theorem \ref{1.33}. Some specific  matrix inversions covered by the $(\alpha,\beta)$-inversion
formula  will be presented in Section 4, among are  three  new matrix inversions related to elliptic divisibility sequence, theta and partial theta functions.

\section{Proof of Lemma \ref{1.3}}
\setcounter{equation}{0}
 Our proof of Lemma \ref{1.3} mainly involves  the following well-known fact about
 the divided difference of polynomials.  See \cite[p.123]{tttsss} for further details.
\begin{yl}\label{yl1} Let $H(x)$ be a polynomial in $x$ of degree no more than
$n-1$ and  $x_0,x_1,\ldots,x_{n}$ be $n+1$ distinct nodes. Then
\begin{align}
 [x_0,x_1,\ldots,x_{n}]H=\sum_{0\leq i\leq n}\frac{H(x_i)}
     { \prod_{j=0,j\neq i}^{n} (x_i-x_j)}=0,
\end{align}
  where  the classical divided  difference of $H(x)$ with respect to
  $\{x_i|0\leq i\leq n\}$
  is recursively defined by
  \begin{align}
  & [x_0]H=H(x_0), \nonumber\\
   & [x_0,x_1]H=\frac{H(x_0)-H(x_1)}{x_0-x_1}, \nonumber\\
    &[x_0,x_1,x_2]H=\frac{[x_0,x_1]H-[x_1,x_2]H}{x_0-x_2},\nonumber\\
&\cdots\qquad\qquad\qquad\cdots\nonumber\\
    &[x_0,x_1,\ldots,x_n]H=\frac{[x_0,x_1,\ldots,x_{n-1}]H-[x_1,
    x_2,\ldots,x_n]H}{x_0-x_n}.
\nonumber\end{align}
\end{yl}
Now write the polynomial $H(x)$  of degree $n-1$ as
$\prod_{i=1}^{n-1}(x+a_i).$ Then Lemma \ref{yl1} is therefore rephrased
explicitly
\begin{align}
 [x_0,x_1,\ldots,x_{n}]H=\sum_{0\leq i\leq n}\frac{\prod_{j=1}^{n-1}(x_i+a_j)}
     { \prod_{j=0,j\neq i}^{n} (x_i-x_j)}=0.\label{diffeq}
\end{align}

Now we are in a good position to show Lemma \ref{1.3}  after Lemma \ref{yl1} given.

\pf  It only needs to check that (\ref{matrixinvdef}) is true for all $n\geq k$. In the case  $n=k$,  it is self-evident.
We only need to consider the case $n>k$.  As such, we compute in a straightforward way
  \begin{align*}
&\sum_{k\leq i\leq n}F(n,i)G(i,k)=\sum_{k\leq i\leq
 n}\frac{b_n}{b_i} \prod_{j=i+1}^{n}m_j\prod_{j=i+1}^{n} \frac{(
s_i-s_{j-1}+a_{j-1}b_{j-1}m_{j-1}) }
      { s_i - s_j  }\\
&\qquad\qquad\qquad\qquad\quad\times\frac{a_{k}}{a_i}\prod_{j=
k}^{i-1}m_j\prod_{j = k}^{i-1}\frac{
      (s_i-s_{j+1}+
        a_{j+1}b_{j+1}m_{j+1})}{
      s_i-s_j}\\
&=b_na_k\prod_{j=k}^{n}m_j\sum_{k\leq i\leq
 n}\frac{1}{a_ib_im_i}\frac{\prod_{j=i}^{n-1} (
s_i-s_{j}+a_jb_jm_j)\prod_{j = k+1}^{i} (s_i-s_{j}+
       a_jb_jm_j) }
     {\prod_{j=k,j\neq i}^{n} (s_i - s_j)}.
\end{align*}
After a bit simplification, we obtain
\begin{align*}
\sum_{k\leq i\leq n}F(n,i)G(i,k)&=b_na_k\prod_{j=k}^{n}m_j\sum_{k\leq i\leq n}\frac{\prod_{j = k+1}^{n-1}
(s_i-s_{j}+
        a_jb_jm_j) }
     { \prod_{j=k,j\neq i}^{n} (s_i - s_j)}\\
&=b_na_k\prod_{j=k}^{n}m_j\sum_{0\leq i\leq n-k}\frac{\prod_{j = 1}^{n-k-1}
(s_{i+k}-s_{j+k}+
        a_{j+k}b_{j+k}m_{j+k}) }
     { \prod_{j=0,j\neq i}^{n-k} (s_{i+k} - s_{j+k})}.
\end{align*}
Observe that
 \begin{align*}
\sum_{0\leq i\leq n-k}\frac{\prod_{j = 1}^{n-k-1}
(s_{i+k}-s_{j+k}+
        a_{j+k}b_{j+k}m_{j+k}) }
     { \prod_{j=0,j\neq i}^{n-k} (s_{i+k} - s_{j+k})}
\end{align*}
is just a special case of Eq.(\ref{diffeq}) under the specifications that $n\to n-k$ and
 $$(x_i,a_i)\to (s_{i+k}, -s_{i+k}+a_{i+k}b_{i+k}m_{i+k}).$$
Hence we obtain
\[\sum_{k\leq i\leq n}F(n,i)G(i,k)=0.\]
 This gives the complete proof of the theorem. \qed

\begin{rk}
As Krattenthaler pointed out, Lemma \ref{1.3} can be derived by use of his  well-known inversion formula. We refer the reader to \cite{kratt-1} for further  detail.
\end{rk}
\section{Proof of Theorem \ref{1.33}}
\setcounter{equation}{0}
In this section, we will show via the use of Lemma \ref{1.3} that (\ref{cond20}) is sufficient but not necessary to   Conjecture \ref{Thmainnewab}, i.e.,  the ($\alpha,\beta$)-inversion formula.

\subsection{Proof of Conjecture \ref{Thmainnewab} under \eqref{cond20}}

For this purpose,  it is convenient to introduce
\begin{dy} Let $\{\alpha_{n,k}\}_{n,k\in \Bbb{Z}}$ and $\{\beta_{n,k}\}_{n,k\in \Bbb{Z}}$ be two arbitrary double
index  sequences over $\mathbb{C}$. We say  $\{\alpha_{n,k}\}_{n,k\in \Bbb{Z}}$ and $\{\beta_{n,k}\}_{n,k\in \Bbb{Z}}$
satisfy \emph{the triple sum identity (TSI)} provided that
 for any integers $n,k,l,m$, it holds
      \begin{align}
\alpha_{n,l}\beta_{m,k}+ \alpha_{n,m}\beta_{k,l}+
 \alpha_{n,k}\beta_{l,m}=0.\label{cond222}
\end{align}
 While, they satisfy \emph{the quintuple sum identity (QSI)} if for any integers $x,y,l,m$, it holds
 \begin{align}
& \alpha_{x,l}\,\alpha_{l,y}\,\beta_{x,l}\,\beta_{m,y}+
  \alpha_{x,l}\,\alpha_{l,x}\,\beta_{m,y}\,\beta_{l,y}\label{quintuple}\\
  &\qquad-
  \alpha_{x,y}\,\alpha_{l,y}\,\beta_{x,l}\,\beta_{m,l}
  - \alpha_{x,y}\,\alpha_{l,m}\,\beta_{x,l}\,\beta_{l,y} -
  \alpha_{x,x}\,\alpha_{l,l}\,\beta_{m,y}\,\beta_{l,y}=0.\nonumber
\end{align}
\end{dy}
Later as we will see,   these two identities  are crucial to Theorem \ref{1.33}.   In the following, we proceed to  show that they are in fact equivalent to each others, although both seem very different in form. This equivalency is based on the following two facts.
The first one is that  TSI (\ref{cond222}) is also equivalent to (\ref{cond3}).
\begin{yl}\label{tri}  $\{\alpha_{n,k}\}_{n,k\in \Bbb{Z}}$ and $\{\beta_{n,k}\}_{n,k\in \Bbb{Z}}$  with
 $\beta_{n,k}=-\beta_{k,n}$ satisfy
  TSI (\ref{cond222}) if and only if
  for any integers $l,x,y$,
  \begin{align}
\alpha_{l,x}\beta_{y,l}+ \alpha_{l,y}\beta_{l,x}+
\alpha_{l,l}\beta_{x,y}=0.\label{cond3}
\end{align}
\end{yl}
\pf To show this lemma, it only needs to  derive
 from (\ref{cond3}) the TSI
\begin{align}
\alpha_{n,k}\beta_{x,y}+
 \alpha_{n,x}\beta_{y,k}+
\alpha_{n,y}\beta_{k,x}=0,\label{cond4}
\end{align}
since \eqref{cond4} is the special case of TSI.
 For this, we first see that  as the special case of  (\ref{cond3}), it holds
\begin{align}
\beta_{y,x}=\frac{\alpha_{l,y}\beta_{l,x}-\alpha_{l,x}\beta_{l,y}}{\alpha_{l,l}}.\label{special}
\end{align}
Next, by substituting (\ref{special}) for each $\beta_{x,y}$ in (\ref{cond4}), we obtain
\begin{align*}
\mbox{LHS of (\ref{cond4})}&=\frac{\alpha_{n,k}}{\alpha_{l,l}}\left\{\alpha_{l,x}\beta_{l,y}-\alpha_{l,y}\beta_{l,x}\right\}\\
&+\frac{\alpha_{n,x}}{\alpha_{l,l}}\left\{\alpha_{l,y}\beta_{l,k}-\alpha_{l,k}\beta_{l,y}\right\}+
\frac{\alpha_{n,y}}{\alpha_{l,l}}\left\{\alpha_{l,k}\beta_{l,x}-\alpha_{l,x}\beta_{l,k}\right\}.
\end{align*}
After a series rearrangement, it reduces to
\begin{align}
\mbox{LHS of (\ref{cond4})}&=
\frac{\beta_{l,k}}{\alpha_{l,l}}\left\{\alpha_{n,x}\alpha_{l,y}-\alpha_{n,y}\alpha_{l,x}\right\}\label{cond4-444}\\
&+\frac{\beta_{l,x}}{\alpha_{l,l}}\left\{\alpha_{n,y}\alpha_{l,k}-\alpha_{n,k}\alpha_{l,y}\right\}+\frac{\beta_{l,y}}{\alpha_{l,l}}\left\{\alpha_{n,k}\alpha_{l,x}-\alpha_{n,x}\alpha_{l,k}\right\}.\nonumber
\end{align}
Observe that the left-hand side of (\ref{cond4}) is independent of $l$. This allows us to set $l=n$, reducing the right-hand side of \eqref{cond4-444} to zero.  The lemma is proved.
\qed

The second fact is that
\begin{yl}\label{qui} Two sequences
 $\{\alpha_{n,k}\}_{n,k\in \Bbb{Z}}$ and $\{\beta_{n,k}\}_{n,k\in \Bbb{Z}}$ with $\beta_{n,k}=-\beta_{k,n}$ satisfy
 QSI (\ref{quintuple}) if and only if they satisfy \eqref{cond3}.
\end{yl}
\pf Now that $\{\alpha_{n,k}\}_{n,k\in \Bbb{Z}}$ and $\{\beta_{n,k}\}_{n,k\in \Bbb{Z}}$ satisfy
 QSI (\ref{quintuple}), in which we may take $y=x$  to get
  \begin{align*}
  \alpha_{x,p}\,\alpha_{l,x}\,\beta_{m,x}\left\{\,\beta_{x,l}+\beta_{l,x}\right\}-\alpha_{x,x}\beta_{x,l}\left\{\alpha_{l,x}\beta_{m,l}
  +\alpha_{l,m}\,\beta_{l,x}-\alpha_{l,l}\,\beta_{m,x}\right\}=0,
\end{align*}
which can further be simplified  to
(\ref{cond3}) by the prior requirement that $\beta_{x,l}=-\beta_{l,x}$ and replacing $m$ with $y$. Conversely, suppose (\ref{cond3}) holds. Then making the parametric replacement
$(x, l)\to (l, x)$ in (\ref{cond3}),
we have
\begin{align} \alpha_{x,x}\beta_{l,y}=\alpha_{x,l}\beta_{x,y}-\alpha_{x,y}\beta_{x,l}.\label{210}
\end{align}
Alternatively, setting $x\to m$ in (\ref{cond3}), we have
\begin{align} \alpha_{l,l}\beta_{m,y}\,=\alpha_{l,m}\beta_{l,y}-\alpha_{l,y}\beta_{l,m}.\label{211}
\end{align}
Upon multiplying (\ref{210}) with (\ref{211}), we find
\begin{align*}
&\qquad\alpha_{x,x}\,\alpha_{l,l}\,\beta_{l,y}\,\beta_{m,y}
=(\alpha_{x,l}\beta_{x,y}-\alpha_{x,y}\beta_{x,l})(\alpha_{l,m}\beta_{l,y}-\alpha_{l,y}\beta_{l,m})\\
&=\alpha_{x,l}\alpha_{l,m}\beta_{x,y}\beta_{l,y}-
\alpha_{x,l}\alpha_{l,y}\beta_{x,y}\beta_{l,m}-\alpha_{x,y}\alpha_{l,m}\beta_{x,l}\beta_{l,y}+\alpha_{x,y}\alpha_{l,y}\beta_{x,l}\beta_{l,m}
.
\end{align*}
Upon substituting these relations into the left-hand side of (\ref{quintuple}), we arrive at
\begin{align*}
 \mbox{LHS of (\ref{quintuple})}
 &=\alpha_{ x,l }\beta _{ m,y } \left(\alpha _{ l,y } \beta _{ x,l
 } +\alpha _{ l,x } \beta _{ l,y }\right)-\alpha_{ x,l }\beta _{ x,y }\left(\alpha _{ l,y
 } \beta _{ m,l }+\alpha _{ l,m }\beta _{ l,y
 } \right)\\
 &=\alpha _{ x,l }(\beta _{ m,y }\alpha _{ l,l } \beta _{
 x,y
 } -\beta_{x,y } \alpha _{l,l} \beta_{m,y })=\alpha _{ x,l }\alpha _{ l,l } (\beta _{ m,y }\beta _{
 x,y
 } -\beta_{x,y } \beta_{m,y })=0.
 \end{align*}
 In the  ante-penultimate equality, we have utilized (\ref{cond3}).
This completes the proof of the lemma. \qed

Summing up,  we have
\begin{xinzhi}\label{pro1} Suppose $\{\beta_{n,k}\}_{n,k\in \Bbb{Z}}$ is anti-symmetric, i.e.,
$\beta_{n,k}=-\beta_{k,n}$. If $\{\alpha_{n,k}\}_{n,k\in \Bbb{Z}}$ and $\{\beta_{n,k}\}_{n,k\in \Bbb{Z}}$ satisfy  TSI (\ref{cond222}), then they  satisfy QSI (\ref{quintuple}). Vice versa.
\end{xinzhi}
\pf From Lemmas \ref{tri} and
\ref{qui}, it is obvious that  TSI (\ref{cond222}) and QSI (\ref{quintuple}) are
equivalent to each others.\qed

Now we are ready to show Conjecture \ref{Thmainnewab} under \eqref{cond20}/(\ref{cond222}), to which  we often refer as TSI (\ref{cond222}). In other word, if (\ref{cond222}) holds true,  then both (\ref{cc12}) and (\ref{dd12}) just form a matrix inversion.

\pf   To show Conjecture \ref{Thmainnewab},  it is enough to reformulate  $F(n,k)$ and $G(n,k)$  given by Conjecture \ref{Thmainnewab}   in the form given by Lemma \ref{1.3}. In other word,  we assume \begin{align*}
&\frac{\prod_{i=k}^{n-1} \alpha_{i,k}}
{\prod_{i=k+1}^{n}\beta_{i,k}}:= \frac{b_n}{b_k} \prod_{i=k+1}^{n} \frac{m_i(
s_k-s_{i-1}+a_{i-1}b_{i-1}\,m_{i-1}) }
      { s_k - s_i  }\qquad\quad\mbox{and}\\
&\frac{\alpha_{k,k}}{\alpha_{n,n}} \frac{\prod_{i=k+1}^{n}
\alpha_{i,n}} {\prod_{i=k}^{n-1}\beta_{i,n}}:=\frac{a_{k}}{a_n}\prod_{i = k}^{n-1}\frac{\,m_i\,
      (s_n-s_{i+1}+
        a_{i+1}\,b_{i+1}\,m_{i+1})}{
      s_n-s_i}.
 \end{align*}
Both can be further restated as
  \begin{align}
\prod_{i=k+1}^n
\frac{b_i}{b_{i-1}}\frac{m_i(s_k-s_{i-1}+ a_{i-1}b_{i-1}\,m_{i-1}
) }
      { s_k - s_i  }= \prod_{i=k+1}^n
\frac{\alpha_{i-1,k}}
{\beta_{i,k}},\label{nnn}
 \\ \prod_{i=k}^{n-1}\frac{a_i}{a_{i+1}}\frac{\,m_i\,
      (s_n-s_{i+1}+
        a_{i+1}b_{i+1}\,m_{i+1})}{
      s_n-s_i}= \prod_{i=k}^{n-1}\frac{\alpha_{i,i}} {\alpha_{i+1,i+1}}\frac{\alpha_{i+1,n}}
      {\beta_{i,n}}.\label{mmm}
 \end{align}
 Now we can easily  deduce from (\ref{nnn})  via induction on $n$ that for integers $i$,
 \begin{align*}
 \frac{b_i}{b_{i-1}}\frac{m_i(s_k-s_{i-1}+a_{i-1}b_{i-1}\,m_{i-1}
) }
      { s_k - s_i  }=\frac{\alpha_{i-1,k}}
{\beta_{i,k}}.
 \end{align*}
Analogously, by using  (\ref{mmm}) and by induction on $k$,  we  obtain
 \begin{align*}
\frac{a_i}{a_{i+1}}\frac{\,m_i\,
      (s_n-s_{i+1}+
        a_{i+1}b_{i+1}\,m_{i+1})}{
      s_n-s_i}=\frac{\alpha_{i,i}} {\alpha_{i+1,i+1}}\frac{\alpha_{i+1,n}}
      {\beta_{i,n}}.
 \end{align*}
Now, for our purpose, we  define  that  for a fixed integer $l$,
\begin{subequations}\label{condis}
\begin{align}
 s_n:=\displaystyle\frac{\alpha_{l,n}}{\beta_{l,n}},~~&~~ m_n:=\frac{\alpha_{n,l}}{\beta_{n,l}},\\
 a_n:=-\displaystyle\frac{\alpha_{n,n}}{\alpha_{n,l}},~~&~~
  b_n:=\frac{\alpha_{l,l}}{\alpha_{n,l}}.
   \end{align}
   \end{subequations}
 Subsequently,
    by substituting (\ref{condis}) into \eqref{nnn} and \eqref{mmm} and then making some simplifications,  we finally  achieve
\begin{align*}  \frac {
-\alpha_{l,k}\beta_{l,i-1}\alpha_{i-1,l}+\alpha
_{l,i-1}\beta_{l,k}\alpha_{i-1,l}-\alpha_{i-1,i-1}\alpha_{l,l}\beta_{l,
k}}{ \left(\alpha_{l,k}\beta_{l,i}-\alpha_{l,i}\beta_{l,k}\right)
\beta_{l,i-1}}=\frac{\alpha_{i-1,k}}
{\beta_{i,k}}
 \end{align*}
 and
\begin{align*}
\frac { -\alpha_{l,n}\beta_{l,i+1}\alpha_{i+1,l}+\alpha
_{l,i+1}\beta_{l,n}\alpha_{i+1,l}-\alpha_{i+1,i+1}\alpha_{l,l}\beta_{l,
n}}{\left(\alpha_{l,n}\beta_{l,i}-\alpha_{l,i}\beta_{l,n}
\right)\beta_{l,i+1}}=\frac{\alpha_{i+1,n}} {\beta_{i,n}},
 \end{align*}
both of which turn out to be, after further simplification,
\begin{align}
 \mathcal{L}(i-1,k;l,i)=
\mathcal{L}(i+1,n;l,i)=0,
\end{align}
where $\mathcal{L}(x,y;l,m)$ denotes the sum on the left-hand side of
(\ref{quintuple}). It is asserted by the known condition of QSI (\ref{quintuple}). As Proposition \ref{pro1} shows, the latter is equivalent to TSI (\ref{cond222}). The conjecture is thus confirmed. \qed

\subsection{Why (\ref{cond20}) is not  necessary to Conjecture \ref{Thmainnewab}}
In order to clarify this point, assuming that  (\ref{cond20}) is true while both (\ref{cc12}) and (\ref{dd12})  compose  a matrix inversion, we now  set up two different ways to calculate
$\beta_{k,n}$ provided that $\{\alpha_{k,n}\}_{n,k\in  {\Bbb Z}}$ and $\{\beta_{k,k+1}\}_{k\in {\Bbb Z}}$ are given.

For this purpose, we start with a special case of Lemma \ref{tri}. At first, by (\ref{cond222}), we may obtain an expression for $\{\beta_{k,n}\}_{k\leq n\in {\Bbb Z}}$ in terms of $\{\beta_{n-1,n}\}_{n\in {\Bbb Z}}$ as below.
\begin{dl}\label{prop3.2}  Suppose that $\{\alpha_{k,n}\}_{n, k\in {\Bbb Z}}$ and $\{\beta_{k,n}\}_{n, k\in {\Bbb Z}}$ satisfy TSI (\ref{cond222}),  $\beta_{n,k}=-\beta_{k,n}$. Then for $k\leq n$, it holds
\begin{align}
\beta_{k,n}=\sum_{i=k+1}^n\left(\frac{\alpha_{i-1,k}}{\alpha_{i-1,i-1}}
\prod_{j=i+1}^n
\frac{\alpha_{j-1,j}}{\alpha_{j-1,j-1}}\right)\beta_{i-1,i}.\label{hhh}
\end{align}
\end{dl}
\pf It suffices to make in (\ref{cond222}) the parametric replacement
\[(n,l,m,k)\to (n-1,k,n,n-1).\] We obtain at once
\begin{align}
\beta_{k,n}=\frac{\alpha_{n-1,n}}{\alpha_{n-1,n-1}}
\beta_{k,n-1}+\frac{\alpha_{n-1,k}}{\alpha_{n-1,n-1}}\beta_{n-1,n}.\label{ggg}
\end{align}
At this stage, we  recognize (\ref{ggg})  as a recursive relation with respect to $\{\beta_{k,k+n}|n\geq 1, k~\mbox{fixed}\}$. By iterating this  recurrence  repeatedly $n-k$ times and then we obtain (\ref{hhh}).
\qed

We think that  the expression (\ref{hhh}), whereas not equivalent to TSI (\ref{cond222}), can be taken as an efficient way to search for possible  $(\alpha,\beta)$-inversions.  The next are two such examples which are obtained as two general solutions to TSI (\ref{cond222}) by making use of (\ref{hhh}).
\begin{tl}\label{gensolution-0}
Let  $\{a_n,b_n, x_n,y_n,t_n\}_{n\in {\Bbb Z}}$ be  arbitrary complex sequences. Suppose that  $\{\beta_{k,n}\}_{n, k\in {\Bbb Z}}$ is subject to $\beta_{n-1,n}=t_n,$ $ \beta_{k,n}=-\beta_{n,k}$. Define
\begin{align}
\alpha_{k,n}:=\prod_{i=1}^nx_i\big/\prod_{i=1}^ky_i.\label{secondrec}
\end{align}
Then
$\{\alpha_{k,n}\}_{n, k\in {\Bbb Z}}$ and $\{\beta_{k,n}\}_{n, k\in {\Bbb Z}}$  satisfy TSI (\ref{cond222}) if and only if
\begin{align}
\beta_{k,n}=\sum_{i=k+1}^{n}t_i\frac{\prod_{j=i+1}^n x_j}{\prod_{j=k+1}^{i-1} x_j}.\label{newnewnew}
\end{align}
\end{tl}
\pf  It is clear that if $\{\alpha_{k,n}\}_{n, k\in {\Bbb Z}}$ and $\{\beta_{k,n}\}_{n, k\in {\Bbb Z}}$  satisfy TSI (\ref{cond222}), then the relation  (\ref{newnewnew}) follows from Theorem \ref{prop3.2} directly. Conversely, suppose (\ref{newnewnew}) is known. We only need to
 check
\begin{align}
\alpha_{n,l}\beta_{m,k}+ \alpha_{n,m}\beta_{k,l}+
 \alpha_{n,k}\beta_{l,m}=0.\label{happy}
\end{align}
Without loss of generality, suppose that $l\geq k\geq m$. In view of the arbitrariness of
$\{t_n\}_{n\in \mathbb{Z}}$,  it  only needs to show the coefficients of $t_{M}$ with $M: l\geq M\geq k$ in the sum on the left-hand side of
 (\ref{happy}) is zero. For this, when  written   in full form, \eqref{happy} becomes
\begin{multline*}
\frac{\prod_{i=1}^lx_i}{\prod_{i=1}^ny_i}\sum_{i=m+1}^{k}t_i\frac{\prod_{j=i+1}^k x_j}{\prod_{j=m+1}^{i-1} x_j}
\\+\frac{\prod_{i=1}^mx_i}{\prod_{i=1}^ny_i}\sum_{i=k+1}^{l}t_i\frac{\prod_{j=i+1}^l x_j}{\prod_{j=k+1}^{i-1} x_j}
\\+
\frac{\prod_{i=1}^kx_i}{\prod_{i=1}^ny_i}\sum_{i=l+1}^{m}t_i\frac{\prod_{j=i+1}^m x_j}{\prod_{j=l+1}^{i-1} x_j}=0.
\end{multline*}  Obviously, the coefficients of $t_{M}$  on the left hand is
\begin{multline*}
\frac{\prod_{i=1}^mx_i}{\prod_{i=1}^ny_i}\frac{\prod_{j=M+1}^l x_j}{\prod^{M-1}_{j=k+1} x_j}-
 \frac{\prod_{i=1}^kx_i}{\prod_{i=1}^ny_i}\frac{\prod_{j=M+1}^m x_j}{\prod_{j=l+1}^{M-1} x_j}\\
 =\frac{\prod_{i=1}^mx_i}{\prod_{i=1}^ny_i}\frac{\prod_{j=M+1}^l x_j}{\prod^{M-1}_{j=k+1} x_j}-
\frac{\prod_{i=1}^kx_i}{\prod_{i=1}^ny_i}\frac{\prod_{j=M}^{l} x_j}{\prod_{j=m+1}^{M} x_j}\\
 =\frac{\prod_{j=1}^lx_j-\prod_{j=1}^lx_j}{\prod_{i=1}^ny_i
 \prod^{M-1}_{j=k+1} x_j\prod_{j=m+1}^{M} x_j}=0.
\end{multline*}
It gives the complete proof of (\ref{happy}).
\qed
\begin{tl}\label{gensolution}Let  $\{a_n,b_n, x_n,y_n,t_n\}_{n\in {\Bbb Z}}$ be  arbitrary complex sequences. Suppose that  $\beta_{n,k}=-\beta_{k,n}$, $\beta_{n-1,n}=a_nb_{n-1}-a_{n-1}b_n,$ and define
\begin{align}
\alpha_{k,n}=x_ka_n+y_kb_n.\label{secondrec}
\end{align}
Then
$\{\alpha_{k,n}\}_{n, k\in {\Bbb Z}}$ and $\{\beta_{k,n}\}_{n, k\in {\Bbb Z}}$  satisfy TSI (\ref{cond222}) if and only if
\begin{align}\beta_{k,n}=a_nb_k-a_kb_n.\end{align}
\end{tl}
\pf It can be verified in a straightforward manner.\qed

In the meantime,  in view of the definition (\ref{matrixinvdef}), it is not hard to establish  another expression for $\beta_{k,n}$
\begin{dl}  $\{\alpha_{k,n}\}_{n,k\in {\Bbb Z}}$ and $\{\beta_{k,n}\}_{n, k\in {\Bbb Z}}$ form an  $(\alpha,\beta)$-inversion if and only if for $k\leq n$, it holds
\begin{align}
\beta_{k,n}=-(f(k,n;k)+f(k,n;n))\big/\sum_{i=k+1}^{n-1}g(k,n;i),\label{secondrec}
\end{align}
where\begin{subequations}
\begin{align}
f(k,n;i)&:=(-1)^{n-i}\prod_{\stackrel{k\leq j_1<j_2\leq n}
{j_1,j_2\neq i}}\beta_{j_1,j_2}\prod_{j=k+1}^{n-1}\alpha_{j,i},\label{defone}\\
 g(k,n;i)&:= (-1)^{n-i}
\prod_{\stackrel{k\leq j_1<j_2\leq n,j_1,j_2\neq i}{j_2-j_1\leq n-k-1}}\beta_{j_1,j_2}\prod_{j=k+1}^{n-1}\alpha_{j,i}.\label{deftwo}
\end{align}\end{subequations}
\end{dl}
\pf According to the definition (\ref{matrixinvdef}), it is clear that $\{\alpha_{k,n}\}_{n,k\in {\Bbb Z}}$ and $\{\beta_{k,n}\}_{n, k\in {\Bbb Z}}$ form an  $(\alpha,\beta)$-inversion if and only if for $n>k$, it holds
\begin{align}
\sum_{i=k}^n\frac{\prod_{j=i}^{n-1} \alpha_{j,i}}
{\prod_{j=i+1}^{n}\beta_{j,i}}~ \frac{\alpha_{k,k}}{\alpha_{i,i}}
\frac{\prod_{j=k+1}^{i} \alpha_{j,i}}
{\prod_{j=k}^{i-1}\beta_{j,i}}=0.\label{inversedef}
\end{align}
The last identity, after simplified by the relation $\beta_{k,n}=-\beta_{n,k}$, is equivalent to
\begin{align}
\sum_{i=k}^nf(k,n;i)=0, \label{secondrec0}
\end{align}
where $f(k,n;i)$ is given by (\ref{defone}).
Further, we split the sum on the left-hand side of (\ref{secondrec0})  into two parts according as the summand $f(k,n;i)$   contains the factor $\beta_{k,n}$ or not. The result is as follows:
 \begin{align*}
f(k,n;k)+\sum_{i=k+1}^{n-1}f(k,n;i)+f(k,n;n)=0
\end{align*}
 where, for $k+1\leq i\leq n-1$, we have
\begin{align*}
f(k,n;i)=g(k,n;i)\beta_{k,n}
\end{align*}
with $g(k,n;i)$ defined by (\ref{deftwo}).
This leads us to (\ref{secondrec}), being thereby  equivalent to \eqref{inversedef}. The theorem is proved.
\qed

Now we are in a good position to explain  why   TSI (\ref{cond222}), i.e., (\ref{cond20}) is not  necessary to Conjecture \ref{Thmainnewab}.  It is because  both (\ref{hhh}) and (\ref{secondrec}) are two recursive relations for $\{\beta_{k,k+n}|n\geq 1\}$. Once $\{\alpha_{k,n}\}_{n,k\in  {\Bbb Z}}$ and $\{\beta_{k,k+1}\}_{k\in {\Bbb Z}}$ are  given as the initial conditions, these two recursive relations may more often than not  produce two different solutions. This  contradicts to the uniqueness of $\beta_{k,n}$ as the inverse of $\big(\alpha_{n,k}\big)_{n\geq k\in  {\Bbb Z}}$.

The following is a short  \emph{Mathematica} program to find $\beta_{k,n}$ recursively via (\ref{ggg}) and (\ref{secondrec}).
\begin{sourcelist}
\begin{verbatim}
d[k,n] denotes the \beta_{k,n} defined by (3.12) while ti[k] for \beta_{k,k+i}
 by (3.14); a[k, n] is \alpha_{k,n}. Note that t1[k],  a[k, n] are all initial
conditions.
---------------------------------------
c[i_]:=t1[i-1]
d[k_,n_]:=Sum[a[i-1,k]/a[i-1,i-1]*Product[a[j-1,j]/a[j-1,j-1],{j,i+1,n}]*c[i],
          {i,k+1,n}]
t2[k_]:=(a[1+k,2+k]t1[k]+a[1+k,k]t1[1+k])/a[1+k,1+k]
t3[k_]:=(a[1+k,3+k]a[2+k,3+k]t1[k]t2[k]t1[1+k]-a[1+k,k]a[2+k,k]t1[1+k]
      t2[1+k]t1[2+k])/(a[1+k,2+k]a[2+k,2+k]t1[k]t2[1+k]-a[1+k,1+k]a[2+k,1+k]
      *t2[k]t1[2+k])
t4[k_]:=(a[1+k,4+k]a[2+k,4+k]a[3+k,4+k]t1[k]t2[k]t3[k]t1[1+k]t2[1+k] t1[2+k]
      +a[1+k,k] a[2+k,k]a[3+k,k]t1[1+k]t2[1+k]t3[1+k]t1[2+k] t2[2+k]
      *t1[3+k])/(a[1+k,3+k]a[2+k,3+k] a[3+k,3+k]t1[k]t2[k]t1[1+k] t3[1+k]
      *t2[2+k]-a[1+k,2+k]a[2+k,2+k]a[3+k,2+k]t1[k]t3[k]t2[1+k]t3[1+k]
      *t1[3+k]+a[1+k,1+k]a[2+k,1+k]a[3+k,1+k]*t2[k]t3[k]t1[2+k] t2[2+k]t1[3+k])
\end{verbatim}
\end{sourcelist}

As an example, we list some computational results to justify our argument.

\begin{lz} Set $\alpha_{k,n}=k+n$ and $\beta_{k,k+1}=k$. Then the output by the above program
are
\begin{align*}
t2[k]-d[k,k+2]&=0,\\
t3[k]-d[k,k+3]&=\frac{8 k^3+32 k^2+32 k+5}{8 k^3+36 k^2+52 k+24}, \\
t4[k]-d[k,k+4]&=\frac{2 k+7 }{8 (k+1) (k+2) (k+3) (2 k+3) (2 k+5)}\frac{f(k)}{g(k)},
\end{align*}
where
\begin{align*}
&f(k)=3072 k^{11}+56320 k^{10}+451904 k^9+2085376 k^8+6115168 k^7+11884320 k^6\\
&\qquad+15498308 k^5+13457624 k^4+7592100 k^3+2669648 k^2+540883 k+47328\\
&g(k)=48 k^7+544
   k^6+2452 k^5+5656 k^4+7216 k^3+5232 k^2+2175 k+464.
 \end{align*}
\end{lz}

\section{Some explicit matrix inversions}
\setcounter{equation}{0}
To justify possibly applications of the ($\alpha,\beta$)-inversion in Conjecture \ref{Thmainnewab}, we now list  some important
 concrete inversions via the use of Corollaries \ref{gensolution-0} and \ref{gensolution}.

There comes first is Gasper's matrix inversion
 which appeared in the bibasic
 hypergeometric series.   Gasper obtained such a pair of matrix inversion in
 his extension of Euler's transformation formula. Displayed as  below, it  is indeed a special case of the $(\alpha,\beta)$-inversion formula.
 \begin{lz}[Cf. {\rm  \cite[Eqs.(3.1)/(3.2)]{8}}]
Let $F=(F(n,k))_{n\geq k\in {\Bbb Z}}$ and $G=(G(n,k))_{n\geq k\in {\Bbb Z}}$ be two matrices with entries given by
\begin{subequations}
\begin{align}
F(n,k)&=(-1)^{n-k} p^{(k-n)k}
\frac{(ap^kq^k,bp^{-k}q^k;q)_{n-k}}{(p,bp^{-n-k}/a;p)_{n-k}}\label{cc123}\qquad\mbox{and}\\
 G(n,k)&=p^{\binom{k}{2}-\binom{n}{2}}\frac{ (1-ap^kq^k)(1-bp^{-k}q^k)(ap^nq^k,bq^kp^{-n};q)_{n-k}
}{(1-ap^nq^k)(1-bp^{-n}q^k)(p,bp^{1-2n}/a;p)_{n-k}},\label{dd123}
\end{align}
\end{subequations}
respectively.
 Then $F=(F(n,k))_{n\geq k\in {\Bbb Z}}$ and $G=(G(n,k))_{n\geq k\in {\Bbb Z}}$ is a
matrix inversion.
 \end{lz}
\pf\,\,\,It only needs to take
\begin{align*}
  \alpha_{i,k} &=  (1-ap^kq^i)(1-bp^{-k}q^i)\qquad{and} \\
  \beta_{i,k} &= (p^i-p^k)(1-bp^{-k-i}/a)
\end{align*}
in the $(\alpha,\beta)$-inversion formula. Clearly, $\beta_{i,k}=-\beta_{k,i}$.
As such, it remains to  check (\ref{cond20}). The related verification is left to the reader.
 \qed

 Another important  $(\alpha,\beta)$-inversion formula is the following result due to
 Schlosser, who has used it   successfully to set
  up transformation formulas of bilateral hypergeometric series.
\begin{lz} [Cf. {\rm \cite[Eqs.(7.18)/(7.19)]{00211}}]
Let $F=(F(n,k))_{n\geq k\in {\Bbb Z}}$ and $G=(G(n,k))_{n\geq k\in {\Bbb Z}}$ be two matrices with entries given, respectively, by
\begin{subequations}
\begin{align}
F(n,k)=\displaystyle\bigg(1/b,\frac{(a+bq^k)q^k}{c-a(a+bq^k)};q\bigg)_{n-k}\bigg/\bigg(q,\frac{(a+bq^k)bq^{k+1}}{c-a(a+bq^k)};q\bigg)_{n-k}
\end{align}
and
\begin{align}
G(n,k)=\displaystyle\,\chi_{n,k}
\bigg(q^{k-n+1}/b,\frac{(a+bq^n)q^{k+1}}{c-a(a+bq^n)};q\bigg)_{n-k}\bigg/\bigg(q,\frac{(a+bq^n)bq^{k}}{c-a(a+bq^n)};q\bigg)_{n-k},
\end{align}
where the factor
\[\chi_{n,k}=(-1)^{n-k}q^{\binom{n-k}{2}}\frac{c-(a+bq^k)(a+q^k)}{c-(a+bq^n)(a+q^n)}.\]
\end{subequations}
 Then $F=(F(n,k))_{n\geq k\in {\Bbb Z}}$ and
$G=(G(n,k))_{n\geq k\in {\Bbb Z}}$ is a matrix inversion.
 \end{lz}
 \pf It follows from the $(\alpha,\beta)$-inversion
formula by specifying
\begin{align*}
  \alpha_{i,k} &=  (q^{k}-q^{i}/b)(c-(a+bq^k)(a+q^i)),\\
  \beta_{i,k} &= (q^k-q^{i})(c-(a+bq^k)(a+bq^i)).
\end{align*}
The verification of (\ref{cond20}) is left to the interested reader.
 \qed

One of the most important   matrix inversions  to  elliptic hypergeometric series is  Warnaar's  elliptic
matrix inversion  \cite{1000}.
\begin{lz}
{\rm (Warnaar's elliptic matrix inversion)} Let
$F=(F(n,k))_{n\geq k\in {\Bbb Z}}$ and $G=(G(n,k))_{n\geq k\in {\Bbb Z}}$ be two
infinite lower-triangular with the entries  given by
\begin{subequations}
\begin{align}
F(n,k)&=\frac{\prod_{i=k}^{n-1}
\theta(x_ib_k,x_i/b_k;q)}{\prod_{i=k+1}^n
\theta(b_ib_k,b_i/b_k;q)} \qquad\mbox{and}\\
G(n,k)&=\frac{b_k\theta(x_kb_k,x_k/b_k;q)}
{b_n\theta(x_nb_n,x_n/b_n;q)}\frac{\prod_{i=k+1}^n\theta(x_ib_n,x_i/b_n;q)}
{\prod_{i=k}^{n-1}\theta(b_ib_n,b_i/b_n;q)}.
\end{align}
\end{subequations}
Then $F=(F(n,k))_{n\geq k\in {\Bbb Z}}$ and $G=(G(n,k))_{n\geq k\in {\Bbb Z}}$ is a matrix inversion.
\end{lz}
\pf It suffices to specify in the $(\alpha,\beta)$-inversion
formula
\begin{align*}
\alpha_{i,k}=b_k\theta(x_ib_k,x_i/b_k;q\big)\,\,\mbox{and}\,\,
\beta_{i,k}=b_k\theta(b_ib_k,b_i/b_k;q\big).
\end{align*}
 It is easy to check that
 $\beta_{i,k}=-\beta_{k,i}$ and TSI (\ref{cond20}) is in agreement  with
 the well-known Weierstrass theta identity \cite[Ex. 2.16(i)]{10}:
\begin{align*}
\theta(xy,x/y,uv,u/v;q)-\theta(xv,x/v,yu,u/y;q)=\frac{u}{y}\theta(xu,x/u,yv,y/v;q).
\end{align*}
The conclusion is proved.
\qed

It is worth mentioning  that the above three pairs of matrix inversions are also special $(f,g)$-inversion. The following  three new matrix inversions  are   not special cases of the $(f,g)$-inversion but covered by the $(\alpha,\beta)$-inversion. This fact shows that  the $(\alpha,\beta)$-inversion is essentially different from the $(f,g)$-inversion.
As a matter of fact, by Corollary \ref{gensolution-0}, we may obtain the first new
elliptic
matrix inversion.
\begin{tl} Let
$F=(F(n,k))_{n\geq k\in {\Bbb Z}}$ and $G=(G(n,k))_{n\geq k\in {\Bbb Z}}$ be two
infinite lower-triangular with the entries  given by
\begin{subequations}
\begin{align}
F(n,k)&=\prod_{i=k+1}^n\frac{1}{
(x;q,p)_i(y;q,p)_{i-1}S_{i,k}} \qquad\mbox{and}\\
G(n,k)&=\prod_{i=k}^{n-1}\frac{1}
{(x;q,p)_{i+1}(y;q,p)_{i}S_{i,n}}.
\end{align}
\end{subequations}
Then $F=(F(n,k))_{n\geq k\in {\Bbb Z}}$ and $G=(G(n,k))_{n\geq k\in {\Bbb Z}}$ is a matrix inversion. Here, for any  complex numbers $x,p,q$,  we define
\begin{align}
 S_{k,n}=\sum_{i=k+1}^{n}\frac{t_i}{(x;q,p)_{i-1}(x;q,p)_i}.
   \end{align}\end{tl}
\pf It suffices to take in Corollary \ref{gensolution-0} that
\[x_i= \theta(xq^{i-1};p),\,\,y_i=\theta(yq^{i-1};p).\]  Note that
\begin{align*}
\alpha_{i,k}=\frac{(x;q,p)_k}{(y;q,p)_i},\,\,  \beta_{i,k}=(x;q,p)_i(x;q,p)_kS_{i,k}.
\end{align*}
It leads to the desired inversion.
\qed

In the same line,  we can find the second new  matrix inversion
arising from Schilling and Warnaar's partial theta function  identity.
\begin{tl} Let
$F=(F(n,k))_{n\geq k\in {\Bbb Z}}$ and $G=(G(n,k))_{n\geq k\in {\Bbb Z}}$ be two
infinite lower-triangular with the entries  given by
\begin{subequations}
\begin{align}
F(n,k)&=\frac{\prod_{i=k}^{n-1}(a_i+\Theta(b_k))}{\prod_{i=k+1}^n
(b_i-b_k)L(b_i,b_k)} \qquad\mbox{and}\\
G(n,k)&=\frac{a_k+\Theta(b_k)}{a_n+\Theta(b_n)}
\frac{\prod_{i=k+1}^n(a_i+\Theta(b_n))}
{\prod_{i=k}^{n-1}(b_i-b_n)L(b_i,b_n)}.
\end{align}
\end{subequations}
Then $F=(F(n,k))_{n\geq k\in {\Bbb Z}}$ and $G=(G(n,k))_{n\geq k\in {\Bbb Z}}$ is a matrix inversion. Here, for any two  complex numbers $x,y$, we define
 \begin{align*}
 L(x,y)&:=-\poq{q,xq,yq}{\infty}
 \sum_{n=0}^\infty\frac{\poq{xy}{2n}}{\poq{q,xq,yq,xy}{n}}q^{n},
 \\
 \Theta(x)&:=\sum_{n\geq 0}(-1)^nq^{n(n-1)/2}x^n.
   \end{align*}\end{tl}
\pf The conclusion follows from the $(\alpha,\beta)$-inversion
formula by specifying
\begin{align*}
   \alpha_{i,k}=a_i+\Theta(b_k)\quad
   \mbox{and}\quad
    \beta_{i,k}=(b_i-b_k)L(b_i,b_k).
\end{align*}
In this case, the validity  of (\ref{cond20}) results from the following partial theta function identity
 \begin{align}
 L(x,y)=\frac{\Theta(x)-\Theta(y)}{x-y},\end{align}
which is Lemma 4.3 of \cite{warnaar00} given by  Schilling and Warnaar.
 \qed

We end this paper by the third new matrix inversion related to the elliptic divisibility sequence $\{W_n\}_{n\in\mathbb{Z}}$ first introduced by M. Ward \cite{ward}. Recall that the elliptic divisibility sequence $\{W_n\}_{n\in\mathbb{Z}}$ is  defined recursively by
\begin{align}\label{EDS}\left\{
                   \begin{array}{ll}
     & W_{n+2}W_{n-2}= W_{n+1}W_{n-1}W_2^2 - W_{1}W_{3}W_n^2 \\
\\
& W_{-n}=-W_n;~ W_0 = 0, W_1 = 1.
                   \end{array}
                 \right.
\end{align}
See \cite{poorten} for details. Our purpose here is to give  a general reciprocal relation for such kind of sequences. It
 convictively shows that the $(\alpha,\beta)$-inversion formula of Conjecture \ref{Thmainnewab} has an advantage over the $(f,g)$-inversion of Theorem \ref{fg} as far as discrete sequences are concerned.
\begin{tl} Let
$F=(F(n,k))_{n\geq k\in {\Bbb Z}}$ and $G=(G(n,k))_{n\geq k\in {\Bbb Z}}$ be two
infinite lower-triangular with the entries  given by
\begin{subequations}
\begin{align}
F(n,k)&=\displaystyle \frac{W_{k}^{2(n-k)}}
{\prod_{i=2k+1}^{n+k}W_{i}\prod_{i=1}^{n-k}W_{i}}\qquad\mbox{and}\\
 G(n,k)&=\displaystyle (-1)^{n-k}\frac{W_k^2}{W_n^2}
\frac{W_n^{2(n-k)}\prod_{i=1}^{n+k-1}W_{i}}
{\prod_{i=1}^{2n-1}W_{i}\prod_{i=1}^{n-k}W_{i}} .
\end{align}
\end{subequations}
Then $F=(F(n,k))_{n\geq k\in {\Bbb Z}}$ and $G=(G(n,k))_{n\geq k\in {\Bbb Z}}$ is a matrix inversion.
\end{tl}
\pf It suffices to take
 in (\ref{quintuplenew}) of Conjecture \ref{Thmainnewab},  preforming as before,
 \begin{align*}
\alpha_{i,j}=W_j^2,\quad \beta_{i,j}=W_{i+j}W_{i-j}.
\end{align*}
Observe that  $\beta_{i,k}=-\beta_{k,i}$, because $W_{-j}=-W_j$ and TSI (\ref{cond20})  is just agreement with
 the property (cf. \cite{poorten}) that for all integers $i,j,k\in {\Bbb Z}$, it holds
\begin{align}\label{property}
   W_i^2  W_{j+k}W_{j-k}+W_j^2  W_{k+i}W_{k-i}+W_k^2  W_{i+j}W_{i-j}=0.
\end{align}
As claimed.
\qed




\end{document}